%Paper: alg-geom/9302001
%From: ANDREATTA%ITNCISCA.bitnet@ICINECA.CINECA.IT
%Date: 01 Feb 1993 13:48:52 +0000

%%%%%%%%%%%%%%%%%%%%%%%%%%%%%%%%
%This is written using Plain TeX. Date:January 28, 1993.
%%%%%%%%%%%%%%%%%%%%%%%%%%%%%%%%

\magnification=1200

% first we put some definitions, mostly for our convenience
\def\proof{\noindent{\bf Proof. }}
\def\remark{\noindent{\bf Remark }}

\def\claim{\noindent{\bf Claim }}

% these are frequently used bold letters and calligraphic
\def\O{{\cal O}}
\def\D{{\cal D}}
\def\E{{\cal E}}
\def\P{{\bf P}}
\def\Q{{\bf Q}}
\def\P{{\bf P}}

\def\C{{\bf C}}
\def\G{{\Gamma}}
\def\H{{\cal H}}

% now morphisms
\def\f{\varphi}
\def\ra{\rightarrow}
\def\iso{\simeq}

%and projectivisations

% and the title
\rightline{January 28, 1993.}
\vskip 2in

{\centerline{\bf A note on non-vanishing and applications}}\par

\medskip
{\centerline{M. Andreatta$^1$  and J. A. Wi\'sniewski$^2$}}
\vglue.3in
\medskip
1)      Dipartimento di Matematica,Universit\'a di Trento, 38050Povo (TN),
Itali
   a
        \par
e-mail :        Andreatta@itncisca.bitnet
\par
2)      Instyt. Matematyki, Uniw. Warszawski, Banacha 2, 02-097 Warszawa,
Poland

\par
e-mail:         Jarekw@plearn.bitnet
\medskip

\beginsection Introduction.

Let $X\subset \P^N$ be a smooth projective variety of dimension $n$
defined over the field of complex numbers.
By $L$ let us denote the restriction of $\O_{\P^N}(1)$ to $X$, by $K_X$
let us denote the canonical divisor of $X$ (a linear equivalence class
of the sheaf $\Omega^n X$ of holomorphic $n$-forms). A (classical) version
of {\sl adjunction theory} concerns studying of an {\sl adjoint} linear system
$\vert K_X+rL \vert$, for a suitably chosen positive integer $r$.
In particular, a typical problem is to decide
whether the {\sl adjoint divisor} $K_X+rL$ is {\sl nef} or
for which value of $r$ it becomes {\sl ample} or even {\sl very ample}.
If $K_X+rL$
is {\sl nef} one may want to find out whether the linear system
$\vert K_X+rL\vert$ has base points, or for which positive integer $m$
its multiple $\vert m(K_X+rL)\vert$ becomes base point free.
\par
Moreover, in the above situation one may release the assumptions concerning
$L$ and $X$, and ask all the above questions if $L$ is merely ample and $X$
is possibly singular.
Problems concerning adjoint divisors has drawn
a lot of attention of algebraic
geometers, starting from the classical works of Castenuovo and Enriques, [CE],
who considered adjoint linear systems on surfaces. Among the references
on these problems we would like to point out the ones which have inspired
us mostly: [BS], [F3], [KMM] and finally the recent paper [K].
\par

The most interesting case concerns situation when the {\sl adjoint
divisor} $K_X+rL$ is {\sl nef} but not ample. Then, although
it can not be said much about the system $\vert K_X+rL\vert$ itself,
a Kawamata---Shokurov
Contraction Theorem asserts that some of its multiple, $\vert m(K_X+rL)\vert$,
is base point free for $m\gg 0$ and defines an {\sl adjoint
contraction morphism}
$\f: X\ra Z$ onto a normal projective variety $Z$, with $\f_*\O_X=\O_Z$.
Understanding of this map seems to be very important for
any classification theory of higher dimensional manifolds.
\par

In the present paper we study the situation when $L$ is
merely ample and the {\sl adjoint
contraction morphism} has fibers of ``small'' dimension.
This last hypothesis allows us to apply an inductive method
which is typical of this theory,
which was called {\sl ``Apollonius method''} by Fujita in [F3].
In the present paper we call it {\sl horizontal
slicing} argument (sometime it is called simply {\sl slicing} but here we need
to distinguish it from {\sl vertical slicing} we also use);
it can be briefly summarized as follows.
\par

Consider a general divisor $X'$ from the linear system $\vert L\vert$
(a hyperplane section of $X$ if $L$ is very ample) and assume that it
is a ``good'' variety of dimension $n-1$
(i.e. has the same singularities as $X$).
By adjunction, $K_{X'}= (K_X+L)_{\vert X'}$ and by Kodaira-Kawamata-
Viehweg Vanishing
Theorem, if $r>1$, the linear system $\vert m(K_{X'}+(r-1)L) \vert$
is just the restriction of $\vert m(K_X+rL)\vert$, so that
the {\sl adjoint contraction morphism} of $X'$ can be related
to the one of $X$.
Moreover, fibers of the {\sl adjoint morphism} of $X'$ will be usually
of smaller dimension  and an inductive argument can be applied.
\par
The {\sl horizontal slicing} argument requires therefore the existence
of a ``good'' divisor $X'$
in the linear system $\vert L\vert$ (a {\sl rung} in the language of [F3]).
The system, however, for an ample (but not very ample) $L$
may a priori be even empty.
To overcome this first difficulty we introduce a local set-up in which the
base of the {\sl contraction morphism} will be affine. We will benefit
from this situation also because we will be able to choose effective divisors
which are rationally trivial. Then the next point
is to ensure that the divisor $X'$
does not contain the whole fiber in question.
In other words, we want to ensure that base point locus of
$\vert L\vert$ ($L$ may be changed by adding a divisor
trivial on fibers of $\f$) does not contain the fiber.
This is what may be called ``relative non-vanishing''.
(Now we can explain why we use the word ``horizontal'': we are used to think
about the map $\f: X\ra Z$ as going vertically,
then every divisor from an ample linear system cuts every ``vertical''
fiber of $\f$ of dimension $\geq 1$, so it lies ``horizontaly''.)
To prove it we apply a method of J. Koll\'ar, as in [K]; it is an
improved version of the one used by Y. Kawamata in proving the
Base Point Free Theorem.
\par

Moreover if the dimension of fibers of $\f$ is $\leq r$ we can descend with our
induction to the case of {\sl adjoint morphism} $\f$
with 1-dimensional fibers.
In such a case a structural theorem for {\sl adjoint morphisms}
of smooth varieties was obtained by Ando:
his theorem (2.3) from [A1]
asserts, for example, that if $r$ and the dimension
of a fiber of $\f$ is 1 and the map $\f$ is birational then, locally,
the map $\f$ is a blow-down of a smooth divisor to a smooth codimension-2
subvariety of $Z$. The ``relative non-vanishing'' turns out to be
sufficient to prove the extension of the Ando's result in the smooth case.
Namely we prove a structural theorem for contractions
whose fiber is of dimension $r$, see theorem 4.1.
\par

Finally, as it turns out, the Kawamata---Koll\'ar method
together with the horizontal slicing allows
to extend the non-vanishing to ``relative spannedness''.
Namely we prove:

\proclaim Theorem. Let $X$ be a normal variety such that $K_X$
is $\Q$-Cartier and $X$ has at worst log terminal singularities;
let $L$ be an ample line bundle over $X$.
Let $\f:X\ra Z$ be an adjoint contraction
supported by $K_X+rL$. Let $F$ be a fiber of
$\f$. Assume moreover that
$$\matrix{
\hbox{either }\hfill& dimF < r+1\hfill&\hbox{if } dimZ <dimX \hfill\cr
\hbox{ or }\hfill& dimF\leq r+1\hfill&\hbox{if }\f\hbox{ is birational}\hfill
}$$
Then the evaluation morphism
$\f^*\f_*L\ra L$ is surjective at every point of $F$.

\par
As a consequence of this one can apply the {\sl horizontal
slicing (or Apollonius) argument}, in the above hypothesis
on the dimension of the fibers,
on varieties with log terminal singularities.

\par
The paper is organized as follows: In the next section
we recall some definitions
concerning singularities as well as the definition of the adjoint contraction
morphism.
Then, in section 2, we introduce the ``affine set-up'' and we
describe the adjoint contraction locally.
Section 3 is about non-vanishing and its proof. In Section 4 we prove
the extension of Ando's structure theorem for contractions of smooth varieties.
In Section 5 we prove the relative spannedness.
\par

The work on the present paper was done when the second named author visited
the University of Trento in the Spring of 1992 and in the Winter of 1993,
supported by an Italian CNR grant (GNSAGA);
he would like also to acknowledge a support from a Polish grant KBN GR54.
The first named author was partially supported by MURST and GNSAGA.

\par
We both would like to thank Edoardo Ballico for many helpful discussions.

\beginsection 1. Singularities and adjoint morphism.

In the present section we recall some pertinent definitions and
results which will be used in the sequel. Our language is compatible with
this of [Ut], or [K], or also [KMM]

Let $X$ be a projective normal variety of dimension $n$
defined over the field of complex numbers.
Assume that $D=\sum d_i D_i$ is an effective $\Q$-divisor ($D_i$ are
prime Weil divisors) such that $K_X+D$ is $\Q$-Cartier. We say that the
pair $(X,D)$ is {\sl log canonical}
( repectively, {\sl purely log terminal}) if for some
(and therefore any) resolution of singularities $f: Y\ra X$
of the pair $(X,D)$ we can write
$$K_Y=f^*(K_X+D)+\sum e_jE_j$$
where $E_j$'s are irreducible divisors with simple normal crossing
and $e_j > -1$ (resp. $e_j \geq  -1$)
(see e.g. [Ut, Ch. 1] for details).
\par

We say that $X$  has {\sl log terminal singularities} if the pair $(X,0)$ is
pur
   ely log
terminal, where $0$ is the ``zero" divisor; in particular we mean that $X$
is considered with no boundary.

\smallskip

\remark. Log terminal singularities of a normal variety $X$
are rational singularities, therefore they are in particular
Cohen-Macaulay, see [KMM, 0.2.17] and [Ke].

\medskip

Let us recall that a $\Q$-Cartier divisor in $X$
is called {\sl nef} if it has non-negative
intersection with any complete curve in $X$.
Assume that $K_X$ is not nef (that is, it has negative intersection
with a complete curve on $X$).
Let $L$ be an ample line bundle (Cartier divisor) on $X$. If $X$
has at worst log terminal singularities and for a positive
rational number  $r$, $K_X+rL$ is \Q-Cartier and nef
(but possibly not ample) then, by the
Kawamata---Shokurov Base Point Free Theorem,
it is semiample and hence the ring
$$R(X,K_X+rL) := \bigoplus_{m\geq 0} H^0(X,m(K_X+rL))$$
is finitely
generated algebra over ${\C}$. Thus there exists an {\sl adjoint
model} $$Z := Proj(R(X,K_X+rL))$$ and a projective surjective
morphism $\f :X\ra Z$  which is given by sections of
a high multiple of  $K_X+rL$ and such that $\f_*\O_X=\O_Z$
(a contraction with connected fibers). We say that $K_X+rL$ is a good
supporting divisor for the morphism $\f$; in fact, this divisor
is a pull-back of an ample divisor from $Z$. If $dimZ=dimX$ then the map
$\f$ is birational, otherwise we call it of fibre type.
\medskip

Let us discuss a local (affine) structure of the map $\f$:
the projective variety
$Z$ has a natural affine covering by
sets $\D_+(h)\iso Spec(R(X,K_X+rL)_{(h)})$, see [H1, II-2],
where  $R(X,K_X+rL)_{(h)}$ denotes the subring of elements of
degree 0 in the localization of $R(X,K_X+rL)$  with respect
to the ideal generated by a homogeneous form $h$.
Let $Z_h$ denote such an affine variety, by $X_h$ let us denote its
pull-back via $\f$ (via base change). Let $\f_h$, respectively $L_h$,
denote again the restriction of $\f$, respectively $L$, to $X_h$. $X_h$
is then projective over $Z_h$, $L_h$ remains clearly $\f_h$-ample.
$ X_h$ is a variety (open subset of $X$), the sections of
$m(K_X+rL)$ associated to multiples of $h$ do  not vanish on $X_h$ and
therefore $m(K_{X_h}+rL_h)$ is a unit in $PicX_h$. Since $m(K_{X_h}+rL_h)$ is
generated for $m$ sufficiently large (not neccesary divisible) we can
shrink $X_h$ so that actually $K_{X_h}+rL_h$ can be assumed to be isomorphic
to $\O_{X_h}$. The morphism $\f_h:X_h\ra Z_h$  is given by the evaluation
map $H^0(X_h, K_{X_h}+rL)\ra (K_{X_h}+rL_h)$. The map $\f_h$ is then as
follows
$$X_h\ni x \mapsto \hbox{ ideal of sections of }K_{X_h}+rL_h\hbox{
vanishing at }x.$$
Note that this affine definition of $\f$ is
clearly compatible with restricting $\f$ to affine subsets.
\par

Frequently, we will replace $Z_h$ by an affine open subset.
Also let us note that, the vanishing of
higher direct images $R^{i}\f_*{\cal F}$ of a sheaf ${\cal F}$
can be understood as vanishing of $H^{i}(X_h,{\cal F})$.
We will also frequently use the existence of effective (Cartier)
divisors on $X_h$  (subvarieties of codimension 1 given by functions).

\beginsection 2. Affine set-up.

We want to understand adjoint contraction morphisms (adjoint maps) locally.
For this purpose we will ``shrink'' $X$ so that the target of the
map will be affine. This will be our affine set-up. In the present
section we review an affine version of adjunction theory and
check basic properties of the construction. Most of this seems
to be standard but we found no reference for it.
\medskip

Let, for a while, $X$ denote a scheme (over
${\bf C}$ ). Set $\G (X):= Spec(H^0(X,\O_X))$. The evaluation of global
functions yields a morphism $\f_{\G}: X\ra \G(X)$ defined as follows
$$X_h\ni x \mapsto \hbox{ ideal of global functions vanishing at }x.$$
It is clear that definitions of $\G(X)$ and $\f_{\G}$
are functorial.  Also, it is not hard to check the following
\par

\claim (2.1). In the above situation the following holds:
\item{(2.1.1)} if $X$ is affine then $\f_{\G}$ is an isomorphism,
\item{(2.1.2)} $({\f_{\G}})_*\O_X=\O_{\G(X)}$,
\item{(2.1.3)} if $X$ is separable over ${\bf C}$ then $\f$ is separable as
well
    [H1, II.4.6],
\item{(2.1.4)} if $X$ is normal then also $\G(X)$ is normal.
\item{}
\bigskip

\noindent (2.2) Now we are ready to set-up our assumptions:
\par
\item{(2.2.1)} Let $X$ be a normal
variety (over ${\bf C}$)
with at worst log terminal singularities and assume that $K_X$
is $\Q$-Cartier. Assume that $L$ is a line bundle (Cartier divisor) on
$X$ and $K_X+rL$ is Cartier for some rational number $r$.
\par

\item{(2.2.2)} Assume that $K_X+rL$ is a
unit in $PicX$. This is equivalent to say that $K_X+rL$ has a nowhere
vanishing section.
\par

\item{(2.2.3)} Set $\f:= \f_{\G}: X\ra Z:= \G(X)$.
Assume that $Z$ is of finite type over $\C$ and $L$ to
be $\f$-ample (thus $\f$ is projective) and (otherwise stated
differently) $r\geq 0$ ,  $r > 0$ if $dimZ<dimX$
\item{}

\par

\remark (2.3). In the above situation
$Z$ is a normal affine variety and the fibers of
$\f$ are connected.
Note also that in this situation $Z:=Proj(R(X, K_X+rL))$
is isomorphic to the
affine scheme $Spec(H^0(X, K_X+rL))= \G (X)$.
Also, the output of the localization procedure
described in the preceeeding section
is exactly the same.
\par
If $\f :X\ra Z$ satisfies all assumptions named in (2.2) then we will call it
a {\sl local adjoint contraction morphism} supported by $K_X+rL$.
We merely use this name to distinguish the affine situation we want to deal
with
   .

\proclaim Lemma 2.4. {\rm (Vanishing theorem)} Let $\f:X\ra Z$ be a
local contraction morphism supported by $K_X+rL$. Assume that
$t$ is an integer and $t > -r$ if $dimZ<dimX$ or $t\geq -r$
if $\f$ is birational. Then
$$H^{i}(X,tL)=0 \hbox{ for } i > 0.$$

\proof
Note that $- K_X = rL$ and $L$ is $\f$-ample. To prove
the vanishing use then the theorem 1.2.5 and the remark 1.2.6
in [KMM], and in the
birational case  also the theorem 1.2.7 in [KMM].
\smallskip

We will frequently use the following ``slicing'' arguments
\par

\proclaim Lemma 2.5. {\rm (Vertical slicing)}
Let $\f:X\ra Y$ be a local adjoint contraction morphism
supported by $K_X+rL$, X is a normal projective variety with at worst
log terminal singularities.
Assume that $X''\subset X$ is a non-trivial
divisor defined by a global function $h\in H^0(X,K_X+rL)=H^0(X,\O_X)$.
Then for a general choice of $h$,
singularities of $X''$ are not worse than these of $X$,
and any section of $L$ on $X''$ extends to $X$.

\proof  The statement on singularities of $X''$ is just
a version of Bertini theorem; for the Bertini theorem,
as well as for the concept of ``general'' we refer to [J], pages 66-67.
The extension property is given by
the vanishing 2.4. Namely, the cokernel of the restriction map
$H^0(X,L)\ra H^0(X'',L_{\vert X''})$  is
contained in $H^1(X,L\otimes\O(-X''))=
H^1(X,L)$  which vanishes by 2.4.
\medskip

\proclaim Lemma 2.6. {\rm (Horizontal slicing)}
Assume the situation from the previous lemma.
\item{(2.6.1)} Let $X'$ be a general divisor
from the linear system $\vert L \vert$. Then the singularities of
$X'$ outside of the base point locus of $\vert L\vert$ are not worse than
these of $X$ and
any section of $L$ on $X'$ extends to $X$.
\item{(2.6.2)} If $\f':=\f_{\vert X'}$ then $K_{X'}+(r-1)L'$ ($L'$ denoting
$L_{
   \vert X'}$) is
$\f'$-trivial and $L'$ is $\f'$-ample
\item{(2.6.3)} Let $Z':=\G(X')$. If either $r > 1$ for $\f$ of fiber type
or  $r\geq 1$
for $\f$  birational then the induced map $Z'\ra Z$ is a
closed immersion. Therefore
the map $\f$ restricted to $X'$ has connected fibers.
\item{}
\par

\proof (2.6.1) is similar to (2.5), (2.6.2) is clear.
As for (2.6.3), note that
the map $Z'\ra Z$ is defined by a homorphism
of rings $H^0(X,K_X+rL)\ra H^0(X',K_{X'}+(r-1)L)$ which is just
the restriction morphism with cokernel contained in
$H^1(X,K_X+(r-1)L)$; this vanishes because of (2.4).
Then connectedness follows from [H1, III, 11.3].

\beginsection 3. Nonvanishing.

In the present section we prove the following

\proclaim Theorem 3.1. Let $\f:X\ra Z$ be a local adjoint contraction
supported by $K_X+rL$ (see assumptions (2.2)). Let $F$ be a fiber of
$\f$. Assume moreover that
$$\matrix{
\hbox{either }\hfill& dimF < r+1\hfill&\hbox{if } dimZ <dimX \hfill\cr
\hbox{ or }\hfill& dimF\leq r+1\hfill&\hbox{if }\f\hbox{ is birational}\hfill
} \leqno{(3.1.1)}$$
Then the base point locus of $L$,
$Bs\vert L\vert := supp(coker(\f^*\f_*L\ra L))$,  does not
contain any component of $F$. Equivalently, there exists a section
of $L$ which does not vanish on any component of $F$.

\remark (3.1.1). If $dimZ=0$ then the theorem is easy: one writes down
Hilbert polynomial $\chi(t):= \chi (X,tL)$ and finds out that
it has zeroes at $-dimF+1$,..., $-1$ so that, since $\chi(t)\not\equiv 0$,
by Kodaira vanishing $L$ has a non-zero section.
Essentially, our point is to extend this argument to our situation.
\par
Actually in the case $dimZ=0$ the above method gives more:
a theorem of Kobayashi-Ochiai, assuming just normal singularities on $X$,
says in fact that, if $r= n+1$ then $(X,L) = (\P^n,\O_{\P^n})$,
while if $r = n$ then $(\Q,\O_\Q(1))$, where \Q \ is an hyperquadric
in $\P^{n+1}$.
Moreover if $r > (n-1)$, besides the two above pairs, the only pair $(X,L)$
with $X$ being $\Q$-Gorenstein is the pair $(\P^2,\O(2))$.
For these results see for instance [B-S] or [F4].

\smallskip
\remark (3.1.2). The dimension of the fiber $F$ is bounded from below:
$dimF \geq r-1$ and $dimF\geq \lfloor r \rfloor$ if $\f$ is birational
(where $\lfloor r \rfloor$ is the integral part of $r$);
we obtain these bounds in the course of our proof but they are
also in [F1] or [F4].
\medskip

The proof of the theorem is based on the technique
developed by Kawamata in proving the base point free theorem; more
precisely we proceed as in the recent paper of Koll\'ar; see [K],
section 2.1: {\sl Modified base point freeness method}. Our proof differs
from Koll\'ar's only by a special choice of divisors involved:
\medskip

\claim (3.2). We can choose a divisor $B$ on $X$ which is a pull back of
a divisor on $Z$ such that $(X,B$) is log canonical outside $F$ and
$(X,B)$ is not log canonical at a generic point of any component $F'$ of $F$.
\par
The proof of claim is the same as (2.2.1) from [K], or (18.22) from [Ut],
we sketch it here:
Take some general functions on $Z$ which
vanish on $\f (F)$. Then the pull backs of their divisors contain
$F$. Let $B$ be their sum. By the general choice of these sections
the first part of the claim follows. To see the second part take
an irreducible component $F'\subset F$: blowing up
$F'$ we obtain an exceptional divisor whose discrepancy with the
respect to $K_X+B$ is $< -1$.
\medskip
Let $f: Y\ra X$ be a log resolution of $B$ i.e. $Y$ is smooth and
all relevant divisors are smooth and cross normally. We then write
$$\matrix{
K_Y=f^*K_X+\sum e_iE_i&\hbox{  where, by assumption  }e_i > -1 \hfill \cr
f^*B=\sum b_iE_i$$ &\hbox{  where  }b_i \geq 0 \hbox{ and } \hfill \cr
f^*(\epsilon L)=A+\sum p_iE_i&\hbox{  if  } dimZ< dimX \hbox{, with } \epsilon
>
    0
\hbox{, or}\hfill \cr
f^*\O_X=A+\sum p_iE_i&\hbox{  if }\f\hbox{ is birational}\hfill
}  \leqno(3.3)$$
where in the last two formulas $A$ is assumed to be $\f \circ f$-ample
$\Q$-divi
   sor
and $0\leq p_i \ll 1$.
\par
Let $F'\subset F$ be an irreducible component, define (c.f. [K], (2.1.2))
$$c:=\hbox{ min}\{ {e_i+1-p_i\over b_i}: F'\subset f(E_i), b_i>0\}
\leqno (3.4)$$
By changing the coefficients $p_i$ a little
we can assume that the minimum is achieved for
exactly one index. Let us denote the corresponding divisor by
$E_0$.
\medskip

\claim (3.5) (see [K] Claim 2.1.3).  By the choice of $c$ we have
\item{(i)}  $0 < c < 1$,
\item{(ii)}   $f(E_0) = F'$,
\item{(iii)} if $cb_i - e_i + p_i < 0$ then $E_i$  is
$f$-exceptional,
\item{(iv)} if $cb_i - e_i + p_i \geq 1$ and $i \ne 0$ then $F'$
is not contained in $f(E_i)$.
\item{}
\par
Let us note that the $\Q$-divisor
$$K_Y + A + \sum (cb_i - e_i + p_i)E_i + f^*(tL) $$
is numerically equivalent to $f^*((t-r)L)$ in the birational, and
to $f^*((t-r+\epsilon)L)$ in the fiber case.
We can write
$$\sum (cb_i - e_i + p_i)E_i = E_0 + H'' - H'+ Fr$$
where $E_0$, $H'$, $H''$ are effective divisor
without common irreducible components; $Fr$ is the fractional divisor
with rational coefficients between 0 and 1, 0 along $E_0$, defined by
$Fr = \sum \{cb_i - e_i + p_i \} E_i$
($\{ a \}$ is the fractional part of the real number $a$).
Moreover $H'$ is $f$-exceptional and $F'$ is not contained in $f(H'')$.
\par
Assume that $t$ is an integer such that $t\geq -r$
if $\f$ is birational, or such that  $t\geq -r+\epsilon$ if $\f$
is of fiber type; by construction
$$f^*(tL) - E_0 + H' - H'' - K_Y - Fr$$
is $\f \circ f$-ample. This divisor restricted to $E_0$ is numerically
equivalen
   t to
$$(f^*(tL)+ H' - H'')_{\vert E_0} - K_{E_0} - Fr_{E_0}$$
and is ample.
By the Kawamata-Viehweg-Koll\'ar vanishing theorems we have then
$$H^i(Y, f^*(tL) - E_0 + H' - H'') = 0\hbox{  and  }
H^i(E_0,(f^*(tL) + H' - H'')_{\vert E_0}) = 0 \hbox{  for } i > 0 .
\leqno (3.6)$$
\par
Let us denote $N(t) :=  f^*(tL) + H' - H''$ for brevity.
By the first vanishing the restriction maps
$$H^0(Y,N(t))\ra H^0(E_0,N(t)_{\vert E_0})\leqno (3.7)$$
are surjective for $t\geq -r$ in the birational and for $t > -r$
in the fiber case.
\par
Assume that then there exists a non-zero section $s$ of $N(t)_{\vert E_0}$.
By surjectivity this extends to a non-zero section of $N(t)$ on $Y$;
since $E_0$ is not contained the support of $H''$ we get a section
$s\in H^0(Y,f^*(tL) + H')$ which is not identically zero along $E_0$.
Moreover  $H^0(Y,f^*(tL)+H') = H^0(X,tL)$ since $H'$ is $f$-exceptional,
thus $s$ descends to a section of $tL$ which does not vanish along
$F' = f(E_0)$ (c.f. [K], (2.1.6)).
\par
{}From this we have immediately that $H^0(E_0,N(t)) = 0$ for $-r \geq t < 0$
(respectively $-r < t < 0$ in the fiber case) and if $dimF'>0$.
\par

Let  $\chi(t):=  \chi(E_0, N(t))$ be the Euler-Poincare characteristic; $\chi$
is a polynomial of degree $\leq d = dimf(E_0) = dimF'$.
By (3.6) and the above we have that  $\chi(t) = 0$ for $-1\geq t \geq -r$
(resp. $-1\geq t > -r$), if $\f$ is birational (resp. if $\f$
is of fiber type.).
Therefore $dim F' \geq \lfloor r \rfloor$ (respectively $dimF'\geq (r-1)$).
Since $\chi(0)\geq 0$ and $\chi(t)>0$ for $t \gg 0$
thus since $dimF' \leq r + 1$ (resp. $ < r +1$) it follows that
$H^0(E_0,N(1))\ne 0$.
\par
Therefore by the surjectivity of the
map in (3.7) (and the subsequent paragraph) we obtain a section of
$L$ in a neighborhood of $F$ which does not vanish
along $F'$; this concludes the proof of the theorem.

\beginsection 4. The structure of an adjoint contraction morphism; smooth case.

In the present section we prove a theorem about the local structure
of an adjoint contraction morphism.

Let $\f: X\ra Z$ be a local adjoint contraction morphism supported by
$K_X+rL$, $r>0$. Let $F$ be a fiber of $\f$.
We assume that  $X$ is smooth in a neighbourhood
of $F$. Our result is subsumed
in the following

\proclaim Theorem 4.1.
If $dimF=r$ then $Z$ is smooth at $\f(F)$ and
\item{(i)} if $dimZ=dimX-r$ then $\f$ is a quadric bundle in
the neighbourhood of $F$,
\item{(ii)} if $dimZ=dimX-r+1$ then $r\leq dimX/2$ and $F =\P^r$,
\item{(iii)} if $\f$ is birational then, in a neighbourhood of
$F$, it is a blow-down of a smooth divisor $E\subset X$ to a smooth
subvariety of $Z$.
\item{}
\par

\remark.
We have just seen, (3.1.2), that every non-trivial fiber of $\f$ is of
dimension at least $r-1$ and if $dimF=r-1$ then $dimZ<dimX$.
Fujita proved [F1, (2.12)] that in the latter case
$F$ is irreducible and $F=P^{r-1}$, and $\f$  is a projective
bundle in a neighbourhood of $F$. Therefore, by
semicontinuity of dimensions of
fibers, it follows that if $dimF=r$ then either

\item{---} $dimX-r+1 \geq  dimZ \geq dimX-r$, or

\item{---} $\f$ is birational and all non-trivial fibers of $\f$
in a neighbourhood of $F$
are $r$-dimensional.

\medskip
\remark. The theorem is a generalization of a result of Ando who proved
a version of it for $r=1$, see [A1], Theorem 2.3.
We will rely on his result, the only exception is the case (i)
which was done by Fujita's slicing technique in [ABW].
Cases (ii) and (iii) will be reduced by horizontal slicing
to the situation of $r=1$ when they coincide and are described in
the Ando's theorem. For $n = 3$ the theorem is known; for instance
the fact that for $n=3$ and $r=2$ the case (ii) does not occur was
proved by Sommese in [S], theorem (3.3), in the case $L$ is spanned
for a slightly more general singularities.

\medskip
{}From now on assume that the situation is as in ii)
or iii) of (4.1)
Although Ando formulated his results for elementary contractions
of projective varieties his results hold in our situation, too.
\par

\proclaim Proposition 4.2.  {\rm ( Ando, [A1],  2.3)}
Assume that $\f:X\ra Z$ is a local
adjoint contraction supported by $K_X+L$.
Moreover assume that a fiber $F$ of $\f$ is of dimension 1
and $\f$ is birational.
Then  locally, in a neighbourhood of $F$, $\f$ is a blow down-morphism
of a smooth divisor on $X$ to a smooth codimension 2 subvariety in smooth $Z$.
\par
\proof The argument is the same as in the course of the proof of
theorem 2.3 from [A1], we will only have to ensure that
that the fiber $F$ is irreducible.
Take  a component $C_i$ of the fiber $F$.
By 2.2 and 1.5 from [ibid],
it is isomorphic to $\P^1$ and $K_X.C_i=-1=-L.C_i$.
By deformation theory, see e.g. the proof of 1.1 from [W3],
small deformations of $C_i$ sweep out a divisor in a neighbourhood of $F$,
call the divisor $E_i$.
Intersection of any two such divisors $E_i$ and  $E_j$
is of codimension 2 in $X$ and it has non-empty intersection with $F$.
An intersection $S$ of general $n-2$
divisors from a good supporting linear system $m(K_X+L)$, $m \gg 0$
will be a smooth surface ($X$ can be assumed to be smooth);
the intersection $E_i\cap E_j\cap S$ is non-empty.
On the other hand, both $E_i\cap S$
and $E_j\cap S$ are contracted to points by a map supported by $L+K_S$,
so that $\f_{\vert S}$ is a contraction of $(-1)$-curves.
Therefore $E_i=E_j$ and the fiber $F$ is irreducible.
\par
To conclude the proof we argue like Ando: by [A1], 1.5,
we know that the normal bundle to $F$ is
$\O\oplus\O\dots\O\oplus\O(-1)$,
so $F$ deforms in a smooth family and finally by Nakano criterium
(see [N])
its contraction is a blow-down.
\medskip

We will need a result about the normalization of components of the fiber $F$.

\proclaim Proposition 4.3. {\rm (Fujita [F1], [F2, (2.2)] and Ye,
Zhang [YZ, lemma4]) }
Assume that $X$, $L$, $r$  and $\f$ are as at the beginning of the section.
Let  $F'$ be an irreducible component of a fiber $F$ of $\f$.
\item{(4.3.1)} If $dimF'=r -1$ then the normalization of $F'$ is
isomorphic to $\P^{r-1}$ and the pull-back of $L$ is $\O(1)$.
\item{(4.3.2)} If $dimF'=r > dim(\hbox{general fiber of }\f)$
then the normalization of
$F'$ is isomorphic to $\P^r$ and the pull-back of $L$ is $\O(1)$.
\item{}
\medskip

Now we can proceed with the proof of (2.1). We can assume
$X$ to be smooth, this is preserved by horizontal slicing:

\proclaim Lemma 4.4.  Assume that a section $X'$ of $L$
does not contain any component of the fiber $F$. Then $X'$ is smooth  on $F$.
\par
\proof Pick up a point $x$ on $F\cap X'$. There exists a rational curve $C$
contained in $F$ such that $C.L=1$, $C$ contains $x$ and it is not contained
in $X'$. This follows from the fact that any two points of $F$
can be joined by such a curve C; this in turn follows from the proposition 4.3.
Therefore $C.X'=1$ and thus $X'$ is smooth at $x$.
\par

\proclaim Lemma 4.5. The line bundle $L$ is base point free in a
neigbourhood of $F$, that is $Bs\vert L\vert := supp(coker(\f^*\f_* L\ra L))$
does not meet $F$.
\par
\proof By 4.2 this is true for $r=1$. For $r>1$ this follows by induction
from a horizontal slicing argument, see 2.6.
\par

\proclaim Lemma 4.6. The fiber $F$ is irreducible.
\par
\proof By 4.2 this is true for $r=1$.
Now assume $r=2$ and $F$ contains a component of
dimension 1. Then by base point freeness we can choose $X'\in |L|$
meeting the 1-dimensional component of $F$ not at the intersection
of components. This contradicts 2.6.3 as the map $\f_{\vert X'}$
has a disconnected fiber. Now general case follows by horizontal slicing
argument: i.e. by 4.4 and 2.6.

\medskip

{}From 4.3, 4.5 and 4.6 it follows now
\par
\proclaim Lemma 4.7.  $F$ is isomorphic to $\P^r$ and
$L_{\vert F}$ is isomorphic to $\O(1)$.

\bigskip

Now we discuss the normal bundle of $F$
\proclaim Lemma 4.8.
The normal bundle of $F$ in $X$ is uniform of the splitting
type $(-1,0,....,0)$.
\par
\proof As above we proceed by induction with respect to $r$.
For $r=1$ this follows from Ando's result.
Now assume that the lemma is true for $r-1$. Consider
a divisor $X'$ from $\vert L\vert$ let $F'=F\cap X'$.
We then have the following exact  sequences of vector bundles on any line
$C$ in $F'$
$$0\ra NF'/X'_{\vert C} = \O(-1)\oplus \O\oplus\dots\oplus\O \ra NF'/X_{\vert
C}
    \ra
NX'/X_{\vert C}   = O(1)\ra 0$$
$$0\ra NF'/F_{\vert C} =\O (1) \ra NF'/X_{\vert C}  \ra NF/X_{\vert C}  \ra 0$$
and therefore $NF/X_{\vert C} = \O(-1)\oplus\O\oplus\dots\oplus\O$.
But, since $L$ is spanned on $F$, the hyperplane $F'$ can be
chosen arbitrarily, so that it contains any line $C$.
\par

\proclaim Lemma 4.9. {\rm (Elencwajg, see [W2], Prop. 1.9)}
The bundle $NF/X$ is either decomposable
into a sum of line bundles or, for $r\leq dimX/2$,
it can be isomorphic to $\Omega (1)\oplus \O\oplus\dots\oplus\O$.

\bigskip

\noindent
Conclusion of the discussion of ii).
\par
We have to prove first that $Z$ is smooth.
This is true for $r=1$ by Ando. For $r>1$ we argue by induction:
take a section $X'\in\vert L\vert$. Fibers of $\f$ restricted to $X'$
are connected so the morphism associated to $K_{X'}+(r-1)L$ is again onto $Z$,
which, by induction is now smooth.
\medskip

\proclaim Lemma 4.10. For $r>n/2$ the case ii) from 4.1 does not occur.
\par
\proof Because of the lemma 4.9  $NF/X$ is decomposable.
Consider a hyperplane $F'\subset F$ which is a specialization of
a general fiber. Then $NF'/X$ is a specialization
of a trivial bundle, see e.g. [W1]. This however is impossible
since from the sequence
$$0 \ra NF'/F =\O (1) \ra NF'/X  \ra NF/X =\O(-1)\oplus\O\oplus\dots\oplus\O
\ra 0$$
we see that the second Chern class of $NF'/X$ can not be 0.
\par
(Note that in [S], theorem 3.3, the case $(n,r) = (3,2)$ is excluded;
in particular we can assume that $dimF' > 1$).
\bigskip

\noindent
We finish now the discussion of the case iii).
\par
Let $E$ be the exceptional set of $\f$ .
To conclude, by Nakano criterium [N], we have to prove that $E$
is a smooth divisor in a neighbourhood of $F$, that
it has projective bundle structure and $\O_E(E)_{\vert F}=\O(-1)$.
\par
First note that a small deformation of $F$ has a decomposable normal bundle.
Indeed, take $n-1-r$ general very ample divisors on $Z$ and consider
the intersection of their pull-back to $X$, call the resulting variety
by $X''$. The variety $X''$ is smooth in the neighbourhood of $F$
and the restriction of $\f$ contracts a small deformation $F'$ of $F$,
being now a divisor in $X''$, to a point. By adjunction we find out that
$NF'/X'' =\O(-1)$ so that $NF'/X'=\O(-1)\oplus\O\oplus\dots\oplus\O$.
\par
But now, because of 4.9 it follows that also
$NF/X=\O(-1)\oplus\O\oplus\dots\oplus\O$.
Now we know that $H^1(F,NF/X)=0$, thus the family of deformations of $F$
inside $X$ is parametrized by a variety which is smooth at $F$ and is
of dimension $n-1-r$. Moreover as in the proof of 1.7 in [E]
we see that the incidence variety of projective spaces is immersed into $X$,
the image of it  coincides locally with $E$.
Therefore $E$ is smooth in the neighbourhood of $F$ and has
a projective bundle structure contracted by $\f$.
Finally, by adjunction, we see $\O_E(E)_{\vert F}= \O(-1)$
which concludes our argument.

\bigskip
A consequence of the theorem 4.1 (iii) is the following generalization of
the theorem 1.1 in [ES].
\proclaim Corollary 4.11.
Let $\f: X \ra Y$ be a proper birational morphism between smooth varieties.
Put $S = \{y \in Y | dim\f^{-1}(y) \geq 1 \}$, $E = \f^{-1}(S)$.  Assume that
$(E)_{red}$ is an irreducible divisor and also that all the fibers of $\f_{|E}$
have the same dimension. Then $S$ is smooth and
$\f$ is the blowing up of $Y$ at $S$.
\par

\proof Note first that $b_2(X) - b_2(Y) = 1$;
for this we use the irreducibility of $E$ as well as the smoothness of $Y$.
Let $r$ be the dimension of one fiber of $\f_{|E}$; the corollary follows now
im
   mediately
from the above theorem taking $L = -E$.

\bigskip
\remark (4.12). We can give a better description of the
case (ii) of the theorem 4.1.
Let us define $\E=\f_*(L)$ and let $S$ be the singular locus
of the sheaf $\E$; $S$ is the locus of $r$-dimensional fibers of $\f$.
Outside of $S$ the sheaf $\E$ is locally free of rank $r$.
On the other hand let us
choose a fiber $F$ of dimension $r$; from 4.1 and 4.5 it follows that,
possibly shrinking $Z$ and $X$, we can choose $r+1$ sections of $L$
spanning it on $X$. This yields a morphism of $X\ra \P^{r}_Z=\P^r\times Z$
over $Z$. The image of the morphism is a divisor and moreover
the morphism in an embedding outside $\f^{-1}(S)$. Therefore, by Serre
criterion
of normality and Zariski Main Theorem
this is an embedding. We have the following sequence of sheaves on $\P^{r}_Z$:
$$0\ra \O(1)\otimes \O(-X)\ra \O(1)\ra \O(1)_{\vert X}=L\ra 0$$
and pushing it forward we get a presentation of $\E$
(note that $X\in\vert\O(1)\vert$ so that $\O(1)\otimes \O(-X)=\O$
and the push-forward is right exact):
$$0\ra \O_Z\ra\O_Z^{r+1}\ra\E\ra 0.$$
\par

Let $X':= Proj_Z(Sym(\E))$. The surjection $\O_Z^{r+1}\ra\E\ra 0$
yields an embedding $X'\subset\P^r_Z$
over $Z$; moreover the homogeneous ideal of $X'$ in $\P^r_Z$
is generated by a linear function and thus it is irreducible and reduced, and
consequently $X\iso X'$. The isomorphism can described globally as follows:
since $L$ is relatively spanned,
from the evaluation morphism
$\f^*\f_*L=\f^*\E\ra L$ we have a regular
morphism $\psi : X\ra X'$ such that $\psi^*\O(1) = L$,
where $\O(1) $ is the relative ample sheaf on the projectivization
(see [H1, section II.7]). Moreover, the pair
$(X',\O(1))$ is
isomorphic, via $\psi$, to $(X,L)$.
\par

Let us also note that if $X$ is projective and
$L$ is ample then also $\E$ is,
where ``ample'' is understood
as in [H2], section 2.
Indeed, the proof of this statement
is exactly the same as in [H2]
(it is stated there for locally free sheaves only but the arguments
remain the same).
\medskip

\remark (4.13).
We have a natural geometric construction
which is associated to the situation
described in the theorem 4.1.(ii). Let $S\subset Z$ be the set of fibers of
dime
   nsion $r$;
$S$ is the set of singularities of $\E$. The set $S$ is of codimension
$r+1$ in $Z$, at least.
\par
Consider a component $\H$ of the Hilbert scheme of $X$ parametrizing
deformations of a general fiber of $\f$.
We claim that $\H$ is smooth. Indeed, this follows because
the normal bundle of any deformation of such a fiber is either
trivial or $\Omega\P^{r-1}(1)\oplus\O(1)\oplus\O^{(n-2r+1)}$.
Moreover, we have a birational map $\pi: \H\ra Z$ which is biregular
outside $\pi^{-1}(S)$, the fibers over $S$ are $\P^r$'s.
Thus, by the corollary 4.11, $S$ is smooth and $\H$ is
a blow-up of $Z$ along $S$. Considering the universal family
over $\H$, we get the following diagram:
$$\matrix{
Y&\ra&\H\cr
\downarrow&\hfill&\downarrow\cr
X&\ra&Z}$$
where $Y\ra\H$ is a projective bundle and $Y\ra X$ a blow-up of
$X$ along $\pi^{-1}(S)$.

\beginsection 5. Spannedness.

In the present section we want to extend the non-vanishing
argument to prove the following

\proclaim Theorem (5.1). Let $\f:X\ra Z$ be a local adjoint contraction
supported by $K_X+rL$ (see assumptions (2.2)). Let $F$ be a fiber of
$\f$. Assume moreover that
$$\matrix{
\hbox{either }\hfill& dimF < r+1\hfill&\hbox{if } dimZ <dimX \hfill\cr
\hbox{ or }\hfill& dimF\leq r+1\hfill&\hbox{if }\f\hbox{ is birational.}\hfill
} \leqno{(5.1.1)}$$
Then the evaluation morphism
$\f^*\f_*L\ra L$ is surjective at every point of $F$.

The proof will be divided in some lemmata. First of all, by vertical slicing
argument and (3.1), we may assume that the base points of the line bundle $L$
are stricly contained in the fibre $F$; we can also assume that $dimZ \geq 1$,
since in the case
$dim Z = 0$ the theorem is trivial (see the remark (3.1.1)).
In particular we can
assume that the dimension of the base locus of $L$
is less then $(n-1)$.
\par
Now we will estimate singularities of a general divisor from the
linear system of $L$.

\proclaim Lemma (5.2). A general divisor from the linear system
$\vert L\vert$ is reduced and irreducible, and it has at worst
log terminal singularities outside of the base point locus of $L$.

\proof Let $G$ be a general divisor of $L$; the fact that it
has at worst log terminal singularities outside $Bs|L|$ follows from
a Bertini theorem. Again by the Bertini theorem (more precisely:
theorem 6.3. points 3 and 4 of [J]), and
the fact that $Bs|L|$ has codimension at least two, we see that
$G$ is generically reduced and irreducible. Note that since we can assume
that $dimZ \geq 1$ and $dimX \geq 2$, we have that the rational map
defined by $|L|$ has image of dimension $\geq 2$.
On the other hand $G$ has Cohen Macaulay singularities, since $X$ has
log terminal and thus Cohen-Macaulay singularities, therefore
$G$ is reduced.

\proclaim Lemma (5.3). In the situation of the above theorem
there exists a divisor $G$ from $\vert L\vert$
which does not contain any component of the fiber $F$ and which has
at worst log terminal singularities on $F$.
\par

\proof We argue similarly as in the proof of (3.1).
We can assume that $Bs\vert L\vert$
is not empty because otherwise the lemma is true
because of the Bertini theorem.
\par
Let $G$ be a  general divisor from $|L|$ which does not contain any
component of the fiber $F$. We may assume that points
where $G$ is not log terminal are contained in the base locus
of $L$.
\par
Let $\tilde {\Delta}$ be a Cartier divisor on $X$ which is a pull back
of a general divisor on $Z$ passing through $\f (F)$, that is $\tilde {\Delta}$
contains $F$; we may assume that $\tilde {\Delta}$ is
log terminal outside $F$.
\par
As in the proof of the theorem 3.1, let
$f: Y \ra X$ be a log resolution
(i.e. $Y$ is smooth and all relevant divisors are smooth and cross normally).
As before we write
$$\matrix{
K_Y=f^*K_X+\sum e'_iE_i&\hbox{  where, by assumption  }e'_i> -1 \hfill \cr
f^*G=\sum b_iE_i$$ &\hbox{ with }b_i \geq 0 \hfill \cr
f^*\tilde {\Delta}=\sum \delta _iE_i$$ &\hbox{ with }\delta _i \geq 0
\hbox{ and }\hfill \cr
f^* \epsilon L=A+\sum p_iE_i&\hbox{  if  } dimZ< dimX
\hbox{ with } \epsilon > 0 \hbox{ or}\hfill \cr
f^*\O_X=A+\sum p_iE_i&\hbox{  if }\f\hbox{ is birational}\hfill
}  \leqno(5.3)$$
where in the last two formulas $A$ is assumed to be ample $\Q$-divisor
and $0 \leq p_i \ll 1$. Note that $b_i \geq 0$ and $\delta _i \geq 0$;
$f(E_i)\subset G$ if and only if $b_i>0$, while
$f(E_i)\subset \tilde {\Delta}$ if and only if $\delta _i > 0$.
Let $E_1 = \bar G$ be the strict transform of $G$; in particular
$e'_1 = 0$ and $b_1 = 1$.
\par
By adjunction
$$K_{\bar G}=(K_Y+\bar G)_{\vert \bar G}=
(f^*(K_X+G)+\sum_{i \ne 1} (e'_i-b_i)E_i)_{\vert \bar G}$$
and if the divisor $G$ is not log terminal
then, for at least one $i \ne 1$, we have $e'_i-b_i\leq -1$.
Let $S$ be a component of the locus, $W$, of non log terminal points of $G$:
$$ W=\bigcup \{f(E_i): e'_i-b_i \leq -1 \ \rm {with} \ i \not= 1 \}.$$
Note that $W$ is contained in the base locus of $L$ and therefore
in $F$.

\par
Let $\Delta = \alpha \tilde{\Delta}$, where $\alpha$ is a rational number
with $0 < \alpha \ll 1$. Let $e_i = e'_i - \alpha \delta_i$; if
$\alpha$ is sufficiently small then $e_i > -1$ for all $i$
and, for the $i$ such that
$e'_i -b_i > -1$, we have that also $e_i - b_i > -1$.
On the contrary if $i \not= 1$ and  $e'_i-b_i\leq -1$, then
$\delta_i > 0$ and therefore $e_i-b_i < -1$.
This means that the pair $(X,\Delta)$ is purely log terminal and
that $X \setminus W$ is the largest open set such that
$(X,\Delta + G)$ is log canonical.

\par
Again, as in the proof of 3.1, we apply Koll\'ar's argument; let
$$c:=\hbox{ min}\{ {e_i+1-p_i\over b_i}: S\subset f(E_i), b_i>0\}
\leqno (5.5)$$
By changing the $p_i$ we can assume that the minimum is achieved for
exactly one index (note that this index is not $1$).
Let us denote the corresponding divisor by $E_0$.
We have a version of the claim (3.6):

\claim (5.6).
\item{(i)} $0< c < 1$,
\item{(ii)} $f(E_0) = S$.
\item{(iii)} if $cb_i - e_i + p_i \geq 1$ and $i \ne 0$ then $S
\not\subset f(E_i)$.
\item{(iv)} if $cb_i - e_i + p_i < 0$ then $e_i > 0$, hence $E_i$  is
$f$-exceptional,
\item{}
\par

\proof See the proof of the Claim 2.1.3 in [K].

\medskip
Let us note that the $\Q$-divisor
$$K_Y + A + \sum (cb_i - e_i + p_i)E_i + f^*(tL) $$
is numerically equivalent to
$f^*((t-r+c)L) $
in the birational, and
to $f^*((t-r+c+\epsilon)L)$
in the fiber case.
We can write
$$\sum (cb_i - e_i + p_i)E_i = E_0 + H'' - H'+
\sum \{ cb_i - e_i + p_i \} E_i$$
where $E_0$, $H'$, $H''$ are effective
and without common irreducible components and the fractional divisor
has rational coefficients between 0 and 1 and is 0 along $E_0$;
moreover $H'$ is $f$-exceptional and $S$ is not contained in $f(H'')$.
\par
Note that since $c < 1$ we can choose $\epsilon$ such that
$c + \epsilon < 1$. Therefore, if we assume that $t$ is an integer such that
$t\
   geq -r+1$
(both in the birational or fiber case),
then, by construction, the divisor
$$f^*(tL)-E_0+H'-H''-K_Y-Fr$$
is $\f\circ f$-ample.
\par

Since by the previous theorem we have that $dimS <dim F \leq r + 1$ in
the birational or, respectively $dimS<dimF< r+1$ in the fiber case
we obtain a contradiction arguing in the same way as we proved
non-vanishing in the last part of the proof of the previous theorem.
Namely we can produce a section of $L$ which does not vanish along
$S$ which was supposed in the base locus.
\par

\bigskip

{\bf Proof of the theorem 5.1}
Note first that if we have a map betwen noetherian schemes, $X \ra Y$,
with zero dimensional (not necessarly connected) fibers then
any line bundle on $X$ is relatively spanned ([H1], chap III, (3.7)).
In particular this says that the theorem is true for $dimF = 0$ without
any hypothesis on $r$.
Now we want to apply induction with respect to $dimF$.
Let $X'$ be a general divisor from the linear system $\vert L \vert$,
$\f':=\f_{\vert X'}$ and $L'$ = $L_{\vert X'}$.
By the lemma 5.3 $X'$ does not contain any component of $F$ and
it has at worst log terminal singularities while by the
lemma 2.6 any section of $L$ on $X'$ extends to $X$.
If $r > 1$ or $\geq 1$, respectively in the fiber or birational case,
again by lemma 2.6,  $\f:X'\ra Z'$ is a local adjoint contraction
supported by $K_X+(r-1)L$, with fiber $F'$ such that $dimF' = dimF -1$.
\par
Therefore the theorem is proved by induction.

\beginsection References.

\item{[A1]} Ando, T., On extremal rays of the higher
dimensional varieties, Invent. Math. {\bf 81} (1985), 347---357.

\item{[A2]} ------, On the normal bundle of $\P^1$ in a higher
dimensional projective variety, American J. Math. {\bf 113} (1991), 949---961.

\item{[ABW]} Andreatta,M., Ballico, E., Wi\'sniewski, J., Two theorems on
elemen
   tary contractions,
to appear in Math. Ann..

\item{[BS]} Beltrametti, M., Sommese, A.J., On the adjunction theoretic
classification of polarized varieties, J. reine angew. Math. {\bf 427}
(1992), 157---192.

\item{[CE]} Castelnuovo, G., Enriques, F., Sopra alcune questioni fondamentali
n
   ella teoria
delle superficie algebriche, Annali di Mat. pura ed applicata, {\bf sec.3- vol
v
   i}
(1901), 165-225.

\item{[E]} Ein, L., Varieties with small dual varieties, II,
Duke Math J. {\bf 52} (1985), 895---907.

\item{[ES]} Ein, L., Shepherd-Barron, N., Some special Cremona transformations,
Amer. J. Math. {\bf 111-5} (1989), 783---800.

\item{[F1]} Fujita, T., On polarized manifolds whose adjoint bundle
is not semipositive, in {\sl Algebraic Geometry, Sendai}, Adv. Studies
in Pure Math. {\bf 10}, Kinokuniya--North-Holland 1987, 167---178.

\item{[F2]} ------, Remarks on quasi-polarized varieties, Nagoya Math. J.
{\bf 115} (1989), 105---123.

\item{[F3]} ------, Classification theories of polarized varieties,
London Lect. Notes {\bf 115}, Cambridge Press 1990.

\item{[F4]} ------, On Kodaira energy and adjoint reduction of polarized
manifol
   ds,
Manuscripta Mathematica, {\bf 76}, (1991), 59-84.

\item{[H1]} Harshorne, R., {\sl Algebraic Geometry}, Graduate Text in Math.
{\bf 52} Springer Verlag 1977.

\item{[H2]} ------, Ample vector bundles, Publ. I.H.E.S. {\bf 29}
(1966), 319---350.

\item{[J]} Jouanolou, J-P., {\sl Th\'eor\`emes de Bertini et Applications},
Progress in Math. {\bf 42}, Birk\-h\"auser 1983.

\item{[KMM]} Kawamata, Y., Matsuda, K., Matsuki, K.: {\sl Introduction to the
Minimal Model Program} in {\sl Algebraic Geometry, Sendai}, Adv. Studies
in Pure Math. {\bf 10}, Kinokuniya--North-Holland 1987, 283---360.

\item{[Ke]} Kempf, G., Cohomology and convexity, in
{\sl Toroidal embeddings I}, Springer Lecture Notes
in Math. {\bf 339}, (1973), 41---52.

\item{[K]} Koll\'ar, J., Effective base point freeness, preprint (1992).

\item{[N]} Nakano, S., On the inverse of monoidal transformation, Publ. Res.
Inst. Math. Sci., {\bf 6} (1970-71), 483---502.

\item{[S]} Sommese, A.J., On the adjunction theoretic structure of projective
varieties, Complex Analysis and Algebraic Geometry, Proceedings G\"ottingen,
1985 (ed. H. Grauert), Lecture Notes in Math., {\bf 1194} (1986), 175---213.

\item{[Ut]} Koll\'ar et al. Flips and abundance for algebraic threefolds,
A summer seminar at the University of Utah, 1991 notes.

\item{[W1]} Wi\'sniewski, J., On deformation of nef value, Duke Math. J.
{\bf 64} (1991), 325---332.

\item{[W2]} ------, Fano manifolds and quadric bundles, to appear
in Math.Zeit..

\item{[W3]} ------, On contractions of extremal rays of Fano manifolds,
J. reine angew. Math {\bf 417} (1991), 141---157.
\item{[YZ]} Ye, Y. and Zhang, Q., On ample vector bundles whose adjunction
bundles are not numerically effective, Duke Math. J. {\bf 60} (1990),
671---687.

\end